\documentclass{article}
\usepackage{stywhispers}
\usepackage[named]{algo}
\usepackage{graphicx}
\usepackage{amsmath}
\usepackage{amsfonts}
\usepackage{amssymb}
\usepackage{epsfig}
\usepackage{url}
\usepackage{color}

\definecolor{light}{gray}{.93}

\newcommand{\bm}[1]{\mbox{\boldmath $#1$}}

\newtheorem{theorem}{Theorem}

\title{Alternating Direction Algorithms for Constrained Sparse Regression:
\\ Application to Hyperspectral Unmixing}

\name{Jos\'{e} M. Bioucas-Dias  \hspace{0.5cm} {\rm and} \hspace{0.5cm} M\'ario A. T. Figueiredo}
\address{{\it Instituto de Telecomunica\c{c}\~{o}es}, {\it Instituto Superior T\'ecnico},  Lisboa, {\bf Portugal}}

\ninept

\begin{document}

\maketitle

\begin{abstract}
Convex optimization problems are common  in hyperspectral unmixing. Examples include:
the constrained least squares (CLS) and the fully constrained least squares (FCLS)
problems, which are used to compute the fractional abundances in linear mixtures of 
known spectra; the constrained basis pursuit (CBP) problem, which is used to find sparse
(\emph{i.e.,} with a small number of non-zero terms) linear mixtures of spectra from 
large libraries; the constrained basis pursuit denoising (CBPDN) problem, which is a 
generalization of BP that admits modeling errors. In this paper,  we  introduce
two new algorithms to efficiently solve these optimization problems,
based on the alternating direction method of multipliers, a method from
the augmented Lagrangian family. The algorithms are termed
SUnSAL ({\it sparse unmixing by variable splitting  and augmented Lagrangian})
and C-SUnSAL ({\it constrained SUnSAL}). C-SUnSAL solves the CBP and CBPDN
problems, while SUnSAL solves CLS and FCLS, as well as a more general version thereof,
called {\it constrained sparse regression} (CSR).
C-SUnSAL and SUnSAL are shown to outperform off-the-shelf methods in terms of speed and accuracy.

\end{abstract}

\section{Introduction}
{\it Hyperspectral unmixing} (HU) is a {\it source separation problem}  with applications in
remote sensing, analytical chemistry, and other areas \cite{bioucasOberview12, keshava02, Lopes, ke00}.
Given a set of observed mixed hyperspectral vectors, HU aims at
estimating the number of reference spectra (the {\it endmembers}), their
spectral signatures, and their fractional abundances, usually under the assumption
that the mixing is linear \cite{keshava02, ke00}.

Unlike in a canonical source separation problem, the {\it sources}  in HU ({\it i.e.}, the fractional
abundances of the spectra/materials present in the data) exhibit statistical dependency
\cite{JNBD05}. This characteristic, together with the high dimensionality of
the data, places HU beyond the reach of most standard
source separation algorithms, thus fostering  active research in the field.

Most HU methods can be classified as statistical or geometrical \cite{bioucasOberview12}.
In the (statistical) Bayesian framework, all inference relies on the
posterior probability density of the unknowns, given the observations.
According to Bayes' law, the posterior probability density results from
two factors: the observation model (the likelihood), which formalizes
the assumed data generation model, possibly including random perturbations
such as additive noise; the prior, which may impose natural constraints on the
endmembers ({\em e.g.,} nonnegativity) and on the fractional abundances
({\em e.g.,} belonging to the probability simplex, since they are relative
abundances), as well as model spectral variability \cite{Dobigeon:ieeesp:08,
Moussaoui:nc:08, Dias}.

Geometrical approaches exploit the fact that, under the linear mixing model,
the observed hyperspectral vectors belong to a simplex set whose vertices correspond
to the endmembers. Therefore, finding the endmembers amounts to identifying the
vertices of that simplex  \cite{bioucas09, bioucasOberview12, plaza02, Dias, CHANG06, LM07, zare08}.

Sparse regression is another direction recently explored for HU
\cite{bioucasOberview12, conf:daniel:whispers:10}, which has connections with both the statistical
and the geometrical frameworks. In this approach, the problem is formulated as that of
fitting the observed (mixed) hyperspectral vectors with sparse
(\emph{i.e.,} containing a small number of terms) linear mixtures of
spectral signatures from a large dictionary available
\textit{a priori}. Estimating the endmembers is thus
not necessary in this type of methods. Notice that the sparse
regression problems in this context are not standard, as the unknown
coefficients (the fractional abundances) sum to one (the so-called
{\em abundance sum constraint} -- ASC) and are non-negative
({\em abundance non-negativity  constraint}  -- ANC).
These problems are thus referred to as {\it constrained sparse regression} (CSR).

Several variants of the CSR problem can be used for HU; some examples follow.
In the classical {\it  constrained least squares} (CLS)   the fractional abundances
in a linear mixture of known spectra are estimated by minimizing the total squared
error, under the ANC. The {\it fully constrained least squares} (FCLS)  adds the ASC to
the CLS problem.  Although no sparseness is explicitly encouraged
in CLS and FCLS, under some conditions (namely positivity of the spectra) it can be
shown that the solutions are indeed sparse \cite{BrucksteinElad}. {\it Constrained basis
pursuit} (CBP) is a variant of the well-known {\it basis pursuit} (BP)
criterion \cite{chenDonohoSaunders:95}
under the ANC; as in BP, CBP uses the $\ell_1$ norm to explicitly encourage
sparseness of the fractional abundance vectors. Finally, {\it constrained basis
pursuit denoising} (CBPDN) is a generalization of CBP that admits modeling errors ({\it e.g.},
observation noise).

\subsection{Contribution}
In this paper, we introduce a class of alternating direction algorithms to solve
several CSR problems (namely CLS, FCLS,  CBP, and CBPDN). The proposed  algorithms are
based on the {\it alternating direction method of multipliers} (ADMM)
\cite{Glowinski, Gabay, EcksteinBertsekas}, which decomposes a difficult problem
into a sequence of simpler ones. Since ADMM can be derived as a variable splitting
procedure followed by the adoption of an augmented Lagrangian method to solve the
resulting constrained problem, we term our algorithms as SUnSAL ({\it spectral
unmixing by splitting and augmented Lagrangian}) and C-SUnSAL ({\it constrained SUnSAL}).


The paper is organized as follows. Section~\ref{sec:problem_formulation}
introduces notation and formulates the optimization problems. Section 3
reviews the ADMM and the associated convergence theorem. Section 4 introduces
the SUnSAL and C-SUnSAL algorithms. Section 5 presents experimental results, and
Section 6 ends the paper by presenting a few concluding remarks.

\section{Problem Formulation:  CLS, FCLS, CRS, CBP, CBPDN}
\label{sec:problem_formulation}
Let ${\bf A}\in \mathbb{R}^{k\times n}$ denote a matrix containing the $n$
spectral signatures of the endmembers,  ${\bf x}\in\mathbb{R}^n$ denote the
(unknown) fractional abundance vector, and ${\bf y}\in\mathbb{R}^{k}$  be an (observed)
mixed spectral vector.  In this paper, we assume that ${\bf A}$ is known;
this is the case the CSR approach \cite{conf:daniel:whispers:10}, where
${\bf A}$ is a library with a large number of spectral signatures, thus
usually $n > k$. Matrix ${\bf A}$ can also be the output of an
endmember extraction algorithm, in which  case usually $n \ll k$.
The key advantage of the CSR approach is that it avoids the
estimation of endmembers, quite often a very hard problem.

The CLS, FCLS, and  CSR problems are, respectively,  defined as
\begin{align}
  \label{eq:P_21}
  \mbox{$(\mbox{P}_{\mbox{\scriptsize CLS}})$:}\quad &  \min_{\bf x}\;(1/2)\|{\bf Ax}-{\bf y}\|_2^2\\\
       & \mbox{subject to:}\;\;{\bf x} \geq  {\bf 0};\nonumber\\
       \nonumber\\[-0.2cm]
  \label{eq:P_22}
  \mbox{$(\mbox{P}_{\mbox{\scriptsize FCLS}})$:}\quad &  \min_{\bf x}\;(1/2)\|{\bf Ax}-{\bf y}\|_2^2\\\
       & \mbox{subject to:}\;\;{\bf x} \geq  {\bf 0},  \;\; {\bf 1}^T{\bf x} =
       1;\nonumber\\
   \nonumber\\[-0.2cm]
  \label{eq:P_23}
  \mbox{$(\mbox{P}_{\mbox{\scriptsize CSR}})$:}\quad &  \min_{\bf x}\;(1/2)\|{\bf Ax}-{\bf y}\|_2^2+\lambda\|{\bf x}\|_1\\\
       & \mbox{subject to:}\;\;{\bf x} \geq  {\bf 0},\nonumber
\end{align}
where $\|{\bf x}\|_2$ and  $\|{\bf x}\|_1$ denote the  $\ell_2$ and $\ell_1$ norms of $\bf
x$, respectively,  $\lambda \geq 0$ is a parameter controlling the relative weight
between the $\ell_2$ and $\ell_1$ terms, ${\bf 1}$ denotes a column vector of 1's, and the
inequality ${\bf x} \geq {\bf 0}$ is to be understood in the componentwise sense.
The constraints ${\bf x} \geq  {\bf 0}$ and ${\bf 1}^T{\bf x} = 1$ correspond to the
ANC and ASC, respectively.

The CLS problem corresponds to  $\mbox{P}_{\mbox{\scriptsize CSR}}$ with $\lambda = 0$.
The FCLS problem would  also be  equivalent to  $\mbox{P}_{\mbox{\scriptsize CSR}}$ if the ASC was enforced in
(\ref{eq:P_23}), since the   $\ell_1$ norm  would be constant in the  feasible set. The CBP and CBPDN
problems are also equivalent to particular cases of $\mbox{P}_{\mbox{\scriptsize CSR}}$,
as stated next.

The CBP optimization problem is
\begin{eqnarray}
  \label{eq:CBP}
  \mbox{$(\mbox{P}_{\mbox{\scriptsize CBP}})$:} &&  \min_{\bf x}\;\|\bf x\|_1\\
       && \mbox{subject to:}\;\;{\bf Ax =y},\;\; {\bf x}\geq  {\bf 0}.\nonumber
\end{eqnarray}
Notice that $\mbox{P}_{\mbox{\scriptsize CBP}}$ corresponds to $\mbox{P}_{\mbox{\scriptsize CSR}}$ with $\lambda \rightarrow 0$.
The CBPDN optimization problem  is
\begin{eqnarray}
\label{eq:CBPDN}
  \mbox{$(\mbox{P}_{\mbox{\scriptsize CBPDN}})$:} &&  \min_{\bf x}\;\|\bf x\|_1\\
       && \mbox{subject to:}\;\|{\bf Ax -y\|_2}\leq \delta,\;\; {\bf x}\geq  {\bf 0}.
       \nonumber
\end{eqnarray}
Problem $\mbox{P}_{\mbox{\scriptsize CSR}}$ is equivalent to
$\mbox{P}_{\mbox{\scriptsize CBPDN}}$ in the sense that for any
choice of $\delta$ for which $\mbox{P}_{\mbox{\scriptsize CBPDN}}$ is feasible,
there is a choice of $\lambda$ for which the solutions of the two problems coincide \cite{Roc70}.
Finally, notice that $\mbox{P}_{\mbox{\scriptsize CBP}}$
corresponds to $\mbox{P}_{\mbox{\scriptsize CBPDN}}$ with $\delta = 0$.

As in (\ref{eq:P_23}), we do not enforce the ASC in (\ref{eq:CBP}) and (\ref{eq:CBPDN}), as
this  would convert those optimization problems into feasibility ones, since the
objective function would be constant in the feasible set.


\section{The ADMM }
\label{sec:ADMM}
\vspace{-0.1cm}
Consider an unconstrained problem of the form
\begin{equation}
\min_{{\bf x}\in \mathbb{R}^n} \; f_1({\bf x}) + f_2 ({\bf G\, x}),\label{unconstrained_basic}
\end{equation}
where $f_1:\mathbb{R}^n \rightarrow \bar{\mathbb{R}}$, $f_2:\mathbb{R}^p
\rightarrow \bar{\mathbb{R}}$, and ${\bf G}\in\mathbb{R}^{p\times n}$.
The ADMM \cite{EcksteinBertsekas, Gabay, Glowinski},
the key tool in this paper, is as shown in Fig.~\ref{fig:ADMM}.
The following is a simplified version of a theorem of
Eckstein and Bertsekas stating convergence of ADMM.
\begin{theorem}[\cite{EcksteinBertsekas}]
\label{th:Eckstein}{\sl Let ${\bf G}$ have full column rank and
 $f_1, f_2$ be closed, proper, and convex. Consider arbitrary
 $\mu>0$ and ${\bf u}_0, {\bf d}_0\in \mathbb{R}^p$.
Consider three sequences $\{{\bf x}_k \in \mathbb{R}^{n}, \; k=0,1,...\}$, $\{{\bf u}_k
\in \mathbb{R}^{p}, \; k=0,1,...\}$, and $\{{\bf d}_k \in \mathbb{R}^{p}, \; k=0,1,...\}$
that satisfy
\begin{eqnarray}
 {\bf x}_{k+1} & = & \arg\min_{{\bf x}} f_{1}({\bf x})
 + \frac{\mu}{2} \|{\bf G}{\bf x} \! - \!{\bf u}_k \! -\! {\bf d}_k\|_2^2  \label{eq:admm_1}\\
 {\bf u}_{k+1} & = & \arg\min_{{\bf u}} f_{2}({\bf u})
 + \frac{\mu}{2} \|{\bf G}{\bf x}_{k+1} \! - \! {\bf u} \! - \! {\bf d}_k\|_2^2   \label{eq:admm_2}\\
 {\bf d}_{k+1} & = &  {\bf d}_{k} - ({\bf G\, x}_{k+1} - {\bf u}_{k+1}).
\end{eqnarray}
Then, if (\ref{unconstrained_basic}) has a solution, the sequence $\{{\bf x}_k\}$  converges to it; otherwise, at least one of the sequences $\{{\bf u}_k \}$ or $\{{\bf d}_k\}$ diverges.}
\end{theorem}

\begin{figure}[tb]
\begin{center}
\colorbox{light}{\parbox{0.9\columnwidth}{\footnotesize
\begin{algorithm}{ADMM}{
\label{alg:salsa4}}
Set $k=0$, choose $\mu > 0$, ${\bf u}_0$,  and  ${\bf d}_0$.\rule[-0.1cm]{0cm}{0.2cm}\\
\qrepeat \rule[-0.1cm]{0cm}{0.2cm}\\
   $  {\bf x}_{k+1}  \in  \arg\min_{{\bf x}} f_{1}({\bf x})
 + \frac{\mu}{2} \|{\bf G\, x} - {\bf u}_k - {\bf d}_k\|_2^2$\rule[-0.1cm]{0cm}{0.2cm}\\
  $  {\bf u}_{k+1}  \in  \arg\min_{{\bf u}} f_{2}({\bf u})
 + \frac{\mu}{2} \|{\bf G\, x}_{k+1} - {\bf u} - {\bf d}_k\|_2^2$\rule[-0.1cm]{0cm}{0.2cm}\\
     ${\bf d}_{k+1} \leftarrow {\bf d}_{k} - ({\bf G \, x}_{k+1}  - {\bf u}_{k+1})$\rule[-0.1cm]{0cm}{0.2cm}\\
     $k \leftarrow k+1$\rule[-0.1cm]{0cm}{0.2cm}
\quntil stopping criterion is satisfied.
\end{algorithm} }}\vspace{-0.2cm}
 \caption{The alternating direction method of multipliers (ADMM).}\label{fig:ADMM}\vspace{-0.4cm}
\end{center}
\end{figure}

\section{Application of ADMM}

In this section, we  specialize the ADMM to each of the optimization problems
stated in Section \ref{sec:problem_formulation}.

\subsection{ADMM  CSR: the {SUnSAL} Algorithm}

We start by writing the optimization  $\mbox{P}_{\mbox{\scriptsize CSR}}$ in the equivalent form
\begin{eqnarray}
   \min_{\bf x}\;(1/2)\|{\bf Ax}-{\bf y}\|_2^2 + \lambda\|{\bf x}\|_1+\iota_{\mathbb{R}_+^n}({\bf x}),\label{eq:obj_CSR}
\end{eqnarray}
where $\iota_S$ is the indicator function of  the set $S$ ({\em i.e.}, $\iota_S({\bf
x})=0$  if ${\bf x}\in S$ and $\iota_S({\bf x})=\infty$ if ${\bf x}\notin S$). We now
apply the ADMM using the following translation table:
\begin{align}
      f_1({\bf x}) & \equiv \frac{1}{2}\|{\bf Ax}-{\bf y}\|_2^2 \label{eq:f1_CSR}\\
      f_2({\bf x}) & \equiv \lambda\|{\bf x}\|_1+\iota_{\mathbb{R}_+^n}({\bf x})\label{eq:f2_CSR}\\
      {\bf G}  & \equiv {\bf I}.
\end{align}

With the current setting, step 3   of the ADMM (see Fig.~1) requires solving
a quadratic problem, the solution of which is

\begin{align}
       {\bf x}_{k+1} & \leftarrow {\bf B}^{-1}{\bf w}
       \label{eq:line3_CSR}
\end{align}
where
\begin{align}
      {\bf B} & \equiv {\bf A}^T{\bf A+\mu I}\\
      {\bf w} & \equiv {\bf A}^T{\bf y}+\mu ({\bf u}_k +{\bf d}_k).
\end{align}

Step 4 of the ADMM (Fig.~1) is simply
\begin{equation}
    {\bf u}_{k+1} \leftarrow \arg\min_{\bf u} \;(1/2)\|{\bf u}-\bm{\nu}_k\|_2^2+(\lambda/\mu)\|{\bf u}\|_1 + \iota_{\mathbb{R}_+^n}({\bf u}) \label{eq:soft_pos}
\end{equation}
where $\bm{\nu}_k\equiv {\bf x}_{k+1}-{\bf   d}_k$.
Without the term $\iota_{\mathbb{R}_+^n}$, the solution of (\ref{eq:soft_pos}) would be
the well-known soft threshold \cite{chenDonohoSaunders:95}:
\begin{equation}
{\bf u}_{k+1} \leftarrow \mbox{soft}(\bm{\nu}_k,\lambda/\mu).\label{eq:soft}
\end{equation}
A straightforward reasoning  leads
to the conclusion that the effect of the ANC term $\iota_{\mathbb{R}_+^n}$ is to
project onto the first orthant, thus
 \begin{equation}
   \label{eq:soft_pos_sol}
    {\bf u}_{k+1} \leftarrow \max\{{\bf 0},\mbox{soft}(\bm{\nu}_k,\lambda/\mu)\},
\end{equation}
where the maximum is to be understood in the componentwise sense.

\begin{figure}[tb]
\begin{center}
\colorbox{light}{\parbox{0.9\columnwidth}{\footnotesize
\begin{algorithm}{SUnSAL}{
\label{alg:salsa3}}
Set $k=0$, choose $\mu > 0$, ${\bf u}_0$,  and  ${\bf d}_0$.\rule[-0.1cm]{0cm}{0.2cm}\\
\qrepeat \rule[-0.1cm]{0cm}{0.2cm}\\
   ${\bf w} \leftarrow {\bf A}^T{\bf y}+\mu ({\bf u}_k +{\bf d}_k)$\\
   $ {\bf x}_{k+1} \leftarrow {\bf B}^{-1}{\bf w}$\\
   $\bm{\nu}_k \leftarrow {\bf x}_{k+1}-{\bf   d}_k$\\
  $ {\bf u}_{k+1}  \leftarrow \max\{{\bf 0},\mbox{soft}(\bm{\nu}_k,\lambda/\mu)\}$\\
     ${\bf d}_{k+1} \leftarrow {\bf d}_{k} - ({\bf x}_{k+1}  - {\bf u}_{k+1})$\rule[-0.1cm]{0cm}{0.2cm}\\
     $k \leftarrow k+1$\rule[-0.1cm]{0cm}{0.2cm}
\quntil stopping criterion is satisfied.
\end{algorithm} }}\vspace{-0.2cm}
 \caption{Spectral unmixing  by variable slitting and augmented Lagrangian (SUnSAL).}
 \label{fig:SUnSAL}\vspace{-0.4cm}
\end{center}
\end{figure}

Fig.~\ref{fig:SUnSAL} shows the SUnSAL algorithm, which solves the
CSR problem (\ref{eq:P_23}). SUnSAL is obtained by replacing
lines 3 and 4 of ADMM by (\ref{eq:line3_CSR}) and (\ref{eq:soft_pos_sol}), respectively.

The objective function (\ref{eq:obj_CSR}) is
proper, convex, lower semi-continuous, and coercive, thus it has a non-empty
set of minimizers (see \cite{Roc70}, for definitions of these convex analysis concepts). Functions $f_1$ and $f_2$ in (\ref{eq:f1_CSR}) and (\ref{eq:f2_CSR})
are closed and ${\bf G}\equiv{\bf I}$ is
obviously of full column rank, thus Theorem 1 can be invoked to ensure convergence of
SUnSAL.

Concerning the computational complexity, we  refer that,  in hyperspectral applications,
the rank of matrix $\bf B$  is no larger that the number of bands, often of the
order of a few hundred, thus ${\bf B}^{-1}$ can be easily precomputed.
The complexity of the algorithm per iteration is thus $O(n^2)$, corresponding to the matrix-vector
products.

\subsection{ADMM  CLS and FCLS}

To solve the CLS problem, we simply run SUnSAL with $\lambda = 0$. The  algorithm  to solve FCLS problem is also
very similar to SUnSAL, with a modification in step 4  linked to the ASC. To derive the ADMM algorithm
to solve the   FCLS problem, let us write the optimization (\ref{eq:P_22})   in the equivalent form
\begin{eqnarray}
   \min_{\bf x}\;(1/2)\|{\bf Ax}-{\bf y}\|_2^2 + \iota_{\{ 1\}}({\bf 1}^T{\bf x})+\iota_{\mathbb{R}_+^n}({\bf x}),\label{eq:obj_FCLS}
\end{eqnarray}
where $ \iota_{\{ 1\}}({\bf 1}^T{\bf x})$ enforces the ASC. We now
apply the ADMM using the following translation table:
\begin{align}
      f_1({\bf x}) & \equiv \frac{1}{2}\|{\bf Ax}-{\bf y}\|_2^2 + \iota_{\{ 1\}}({\bf 1}^T{\bf x})\label{eq:f1_FCLS}\\
      f_2({\bf x}) & \equiv \iota_{\mathbb{R}_+^n}({\bf x})\label{eq:f2_FCLS}\\
      {\bf G}  & \equiv {\bf I}.
\end{align}
The resulting ADMM algorithm is similar SUnSAL  with $\lambda = 0$, with one difference:
step 3   of the ADMM (see Fig.~1) requires solving a quadratic problem with linear equality constraint,
the solution of which is
\begin{align}
       {\bf x}_{k+1} & \leftarrow {\bf B}^{-1}{\bf w}-{\bf C}({\bf 1}^T{\bf B}^{-1}{\bf w}-1)\label{eq:line3}
\end{align}
where
\begin{align}
      {\bf B} & \equiv {\bf A}^T{\bf A+\mu I}\\
      {\bf C} & \equiv {\bf B}^{-1}{\bf 1}({\bf 1}^T {\bf B}^{-1}{\bf 1})^{-1}\\
      {\bf w} & \equiv {\bf A}^T{\bf y}+\mu ({\bf u}_k +{\bf d}_k).
\end{align}

\begin{figure}[tb]
\begin{center}
\colorbox{light}{\parbox{0.9\columnwidth}{\footnotesize
\begin{algorithm}{SUnSAL (FCLS version)}{
\label{alg:fcls}}
Set $k=0$, choose $\mu > 0$, ${\bf u}_0$,  and  ${\bf d}_0$.\rule[-0.1cm]{0cm}{0.2cm}\\
\qrepeat \rule[-0.1cm]{0cm}{0.2cm}\\
   ${\bf w} \leftarrow {\bf A}^T{\bf y}+\mu ({\bf u}_k +{\bf d}_k)$\\
   $ {\bf x}_{k+1} \leftarrow {\bf B}^{-1}{\bf w}-{\bf C}({\bf 1}^T{\bf B}^{-1}{\bf w}-1)$\\
   $\bm{\nu}_k \leftarrow {\bf x}_{k+1}-{\bf   d}_k$\\
  $ {\bf u}_{k+1}  \leftarrow \max\{{\bf 0},\bm{\nu}_k)$\\
     ${\bf d}_{k+1} \leftarrow {\bf d}_{k} - ({\bf x}_{k+1}  - {\bf u}_{k+1})$\rule[-0.1cm]{0cm}{0.2cm}\\
     $k \leftarrow k+1$\rule[-0.1cm]{0cm}{0.2cm}
\quntil stopping criterion is satisfied.
\end{algorithm} }}\vspace{-0.2cm}
 \caption{SUnSAL for the FCLS problem.}
 \label{fig:SUnSAL2}\vspace{-0.4cm}
\end{center}
\end{figure}

Fig.~\ref{fig:SUnSAL2} shows the FCLS version of the SUnSAL algorithm, which solves the
FCLS problem (\ref{eq:P_22}). We note that, in any SUnSAL version, the ANC  can be
deactivated trivially.

\subsection{ADMM for CBP and  CBPDN: the C-SUnSAL Algorithm}

Given that the CBP problem corresponds to CBPDN with $\delta = 0$, we address
only the latter. Problem $\mbox{P}_{\mbox{\scriptsize CBPDN}}$ is equivalent to
\begin{eqnarray}
   \min_{\bf x}\;  \|{\bf x}\|_1 + \iota_{B(\bf y,\delta)}({\bf Ax}) +\iota_{\mathbb{R}_+}({\bf x}),\label{eq:obj_CBPDN}
\end{eqnarray}
where $B({\bf y},\delta)= \{{\bf z}:\,\|{\bf z-y}\|_2\leq\delta\}$ is a radius-$\delta$ closed
ball around ${\bf y}$.
To apply the ADMM we use the following definitions:
\begin{align}
      f_1({\bf x}) & = 0 \label{eq:f1_CBPDN}\\
      f_2({\bf u}) & = \iota_{B(\bf y,\delta)}({\bf u}_1)+\lambda\|{\bf u}_2\|_1+\iota_{\mathbb{R}_+^n}({\bf u}_2)
      \label{eq:f2_CBPDN}\\
      {\bf G}  & = \bigl[{\bf A}^T\; {\bf I}\,\bigr]^T.\label{eq:G_CBPDN}
 \end{align}
where ${\bf u}  = \bigl[{\bf u_1}^T \,{\bf u}_2^T\bigr]^T$.
With the above definitions, the solution of line 3 of ADMM (see Fig.~1),
a quadratic problem, is
\begin{align}
       {\bf x}_{k+1} & \leftarrow {\bf B}^{-1}{\bf w},\label{eq:bpdn_x}
\end{align}
where
\begin{align}
      {\bf B} & \equiv {\bf A}^T{\bf A + I}\\
      {\bf w} & \equiv {\bf A}^T({\bf u}_{1,k} +{\bf d}_{1,k})+({\bf u}_{2,k} +{\bf d}_{2,k}).
\end{align}

Because the variables ${\bf u}_1$ and ${\bf u}_2$ are decoupled, line 4 of ADMM (Fig.~1)
consists in solving two separate problems,
\begin{align}
    {\bf u}_{1,k+1} & \in \arg\min_{{\bf u}} \;(1/2)\|{\bf
    u} - \bm{\nu}_{1,k}\|_2^2 +\iota_{B(\bf y,\delta)}({\bf u})\label{eq:ball_proj_bpdn}\\
    {\bf u}_{2,k+1} & \in \arg\min_{{\bf u}} \;(1/2)\|{\bf u}-\bm{\nu}_{2,k}\|_2^2+
    (\lambda/\mu)\|{\bf u}\|_1 +\iota_{\mathbb{R}_+^n}({\bf u})\label{eq:soft_pos_bpdn}
\end{align}
where
\begin{align}
    \bm{\nu}_{1,k} & = {\bf Ax}_{k+1}-{\bf d}_{1,k}\\
    \bm{\nu}_{2,k} & = {\bf x}_{k+1}-{\bf d}_{2,k}.
\end{align}
The solution of (\ref{eq:ball_proj_bpdn}) is the projection onto  the ball $B({\bf
y},\delta)$, given by
\begin{equation}
 {\bf u}_{1,k+1}\!\! \leftarrow \psi_B({\bf y},\delta)(\bm{\nu}_{1k})\equiv\left\{
\begin{array}{ll}
  \!\!\!\bm{\nu}_{1k},& \!\!\!\|\bm{\nu}_{1,k}-{\bf y}\|_2\leq \delta\\
  \!\!\!{\bf y}+\frac{\bm{\nu}_{1,k}-{\bf y}}{\|\bm{\nu}_{1,k}-{\bf y}\|_2}\,\delta , &\!\!\!\|\bm{\nu}_{1,k}-{\bf
  y}\|_2 >
  \delta.
\end{array}
\right.\label{eq:u1kplus1}
\end{equation}
Similarly to (\ref{eq:soft_pos_sol}), the  solution of (\ref{eq:soft_pos_bpdn}) is given
by
\begin{align}
    {\bf u}_{2,k+1} & \leftarrow \max\{{\bf 0},\mbox{soft}(\bm{\nu}_{2,k},\lambda/\mu)\}.
    \label{eq:soft_pos_sol2}
\end{align}

\begin{figure}[tb]
\begin{center}
\colorbox{light}{\parbox{0.9\columnwidth}{\footnotesize
\begin{algorithm}{C-SUnSAL}{
\label{alg:salsa2}}
Set $k\leftarrow 0$, choose $\mu > 0$, ${\bf u}_{1,0}$,  ${\bf d}_{1,0}$, ${\bf u}_{2,0}$,  and  ${\bf d}_{2,0}$.\rule[-0.1cm]{0cm}{0.2cm}\\
\qrepeat \rule[-0.1cm]{0cm}{0.2cm}\\
   ${\bf w} \leftarrow {\bf A}^T({\bf u}_{1,k} +{\bf d}_{1,k})+({\bf u}_{2,k} +{\bf d}_{2,k})$\\
   $ {\bf x}_{k+1}  \leftarrow {\bf B}^{-1}{\bf w}$\\
   $\bm{\nu}_{1,k}  \leftarrow {\bf Ax}_{k+1}-{\bf d}_{1,k}$\\
   $ {\bf u}_{1,k+1} \leftarrow \psi_B({\bf y},\delta)(\bm{\nu}_{1,k})$\\
   $\bm{\nu}_{2,k}  \leftarrow {\bf x}_{k+1}-{\bf d}_{2,k}$\\
  $ {\bf u}_{2,k+1}  \leftarrow \max\{{\bf 0},\mbox{soft}(\bm{\nu}_{2,k},\lambda/\mu)\}$\\
   ${\bf d}_{1,k+1} \leftarrow {\bf d}_{1,k} - ({\bf Ax}_{k+1}  - {\bf u}_{1,k+1})$\rule[-0.1cm]{0cm}{0.2cm}\\
   ${\bf d}_{2,k+1} \leftarrow {\bf d}_{2,k} - ({\bf x}_{k+1}  - {\bf u}_{2,k+1})$\rule[-0.1cm]{0cm}{0.2cm}\\
   $k \leftarrow k+1$\rule[-0.1cm]{0cm}{0.2cm}
\quntil stopping criterion is satisfied.
\end{algorithm} }}\vspace{-0.2cm}
 \caption{Constrained spectal unmixing  by variable slitting and augmented Lagrangian (C-SUnSAL).}\label{fig:CSUnSAL}\vspace{-0.4cm}
\end{center}
\end{figure}

Fig. \ref{fig:CSUnSAL} shows the  C-SUnSAL algorithm for CBPDN, which results
from replacing line 3  of ADMM (Fig.~1) by (\ref{eq:bpdn_x}) and line 4 of
ADMM by (\ref{eq:u1kplus1})--(\ref{eq:soft_pos_sol2}). As mentioned above,
C-SUnSAL can be used to solve the CBP problem simply by setting $\delta = 0$.
As in SUnSAL, the ANC  can be deactivated trivially.

The objective function (\ref{eq:obj_CBPDN}) is proper, convex,
lower semi-continuous, and coercive, thus it has a non-empty
set of minimizers. Functions $f_1$ and $f_2$ in (\ref{eq:f1_CBPDN})
and (\ref{eq:f2_CBPDN}) are closed and ${\bf G}$ in (\ref{eq:G_CBPDN}) is
obviously of full column rank, thus Theorem 1 can be invoked to
ensure convergence of C-SUnSAL. Concerning the computational complexity, the scenario is similar to that of
SUnSAL, thus complexity of C-SUnSAL is $O(n^2)$ per iteration.

At this point, we make reference to the work  \cite{conf:szlam2010split}, which has also addressed
the CSR problem  (\ref{eq:P_23})  aiming at hyperspectral applications.  The  algorithm therein 
proposed, although different from SUnSAL, has some similarities that result from  the strong connections
between the {\em split Bregman } framework adopted in  \cite{conf:szlam2010split}
and the ADMM (for these connections see, {\em e.g.}, \cite{esser09}).

\section{Experiments}
\label{sec:experiments}
We now report experimental results  obtained  with simulated data generated according to
${\bf y} = {\bf Ax+n}$,
where ${\bf n}\in \mathbb{R}^k$  models additive perturbations. In hyperspectral applications, these perturbations are mostly model errors  dominated by low-pass components. For this reason,
we generate the noise by  low-pass filtering  samples of zero-mean  i.i.d.  Gaussian  sequences of random variables. We define the signal-to-noise ratio (SNR) as
$$
   \mbox{SNR} \equiv10\log_{10}\left( \frac{\mathbb{E}[\|{\bf Ax}\|_2^2]}{\mathbb{E}[\|{\bf n}\|_2^2]}\right).
$$
The expectations  in the above definition are approximated with sample means over 10 runs.
The original fractional abundance vectors $\bf x$ are generated in the following way:  given $s$,
the number of non-zero components in $\bf x$, we generate  random samples uniformly in the $(s-1)-$simplex  and distribute randomly these $s$ values among  the components of $\bf x$.
We  considered two libraries ({\em i.e.,} matrices $\bf A$): a $200\times 400$  matrix with
zero-mean  unit variance i.i.d. Gaussian entries and a $224\times 498$ matrix with a  selection of 498 materials (different mineral types) from the USGS library denoted splib06\footnote{http://speclab.cr.usgs.gov/spectral.lib06}.

As far as we know,  there are no special purpose algorithms for solving the CSR, CBP, and CBPDN problems. Of course  these  are canonical convex problems, thus they can be tackled with standard convex  optimization techniques. Namely,  the  CLS, which is a particular case of CSR,
can be solved with the MATLAB function  {\tt lsqnonneg}, which we use as baseline in our comparisons.

Tables 1 and 2  report reconstruction SNR (RSNR), defined as
$$
 \text{RSNR} = 10\log_{10}\left(\frac{\mathbb{E}[\|{\bf x}\|^2_2]}{\mathbb{E}[\|{\bf x}-\widehat{\bf x} \|_2^2]}\right),
$$
where $\widehat{\bf x}$ is the estimated fractional abundance vector, and  execution times,
for the two libraries  referred above.  The {\tt lsqnonneg}  is run with its default options.
SUnSAL and C-SUnSAL  run  200 iterations, which was found to be more than enough to achieve convergence.

\begin{table}
\centering
\caption{RSNR values and execution times for the Gaussian library defined in the text (average over 10 runs).}\label{tab:results}
\vspace{0.2cm}
{\footnotesize \begin{tabular}{ c || c | c | c | c | c  | c}
 & \multicolumn{2}{c|}{ SUnSAL} & \multicolumn{2}{|c|}{ C-SUnSAL} & \multicolumn{2}{c}{\tt lsqnonneg } \\
 \hline
SNR  & RSNR &  time & RSNR &  time & RSNR &  time \\
      (dB)      &   (dB)  &  (sec)&  (dB)   & (sec) &   & (sec)  \\  \hline \hline
20  &   {\bf 10} &  0.12 & 3 &  0.12 & 3  & 31\\ \hline
30  &   {\bf 32} &  0.12 & 27 & 0.12 & 25 & 32\\ \hline
40  &  {\bf  37} &  0.12 & 30 & 0.12 & 27 & 48\\ \hline
50  &  {\bf  48} &  0.12 & 47 & 0.12 & 42 & 57\\ \hline
\end{tabular}}
\vspace{-0.4cm}
\end{table}

\begin{table}
\centering
\caption{RSNR values and execution times for the USGS library  (average over 10 runs).}\label{tab:results2}
\vspace{0.2cm}
{\footnotesize \begin{tabular}{ c || c | c | c | c | c  | c}
 & \multicolumn{2}{c|}{ SUnSAL} & \multicolumn{2}{|c|}{ C-SUnSAL} & \multicolumn{2}{c}{\tt lsqnonneg } \\
 \hline
SNR  & RSNR &  time & RSNR &  time & RSNR &  time \\
      (dB)      &   (dB)  &  (sec)&  (dB)    & (sec) &   & (sec)  \\  \hline \hline
30  &   {\bf 6} &   0.13 & 1.5 &  0.13 & -7  & 22\\ \hline
40  &  {\bf  17} &  0.13 & 12.2 & 0.13 &  10 & 32\\ \hline
50  &  {\bf  23} &  0.13 & 14.5 & 0.13 &  15 & 47\\ \hline
\end{tabular}}
\vspace{-0.4cm}
\end{table}

We highlight  the following  conclusions: ({\bf a})
the proposed algorithms achieve  higher accuracy in about two orders of magnitude shorter  time. This is a critical issue in imaging application where an instance of the problem has to be solved for each pixel; ({\bf b}) the lower accuracy  obtained with the USGS matrix is due to the fact that the spectral signatures are highly correlated  resulting in a much harder problem than with the Gaussian matrix.

\vspace{-0.3cm}

\section{Concluding Remarks}
In this paper, we introduced  new algorithms  to solve
a class of  optimization problems arising in spectral unmixing. The proposed  algorithms are
based on the {\it alternating direction method of multipliers}, which decomposes a difficult problem into a sequence of simpler ones.  We  showed that sufficient conditions for convergence are satisfied. In limited  set of experiments, the proposed algorithms were shown to clearly outperform an off-the-shelf optimization tool. Ongoing  work includes a comprehensive experimental  evaluation of the proposed algorithms.


\begin{thebibliography}{10}
\footnotesize


\bibitem{bioucas09}
J.~Bioucas-Dias,
\newblock "A variable splitting augmented Lagrangian approach to linear spectral unmixing",
\newblock in  {\em First IEEE GRSS Workshop on Hyperspectral Image and Signal Processing-WHISPERS'2009}, Grenoble, France, 2009.

\bibitem{bioucasOberview12}
J.~Bioucas-Dias, A. Plaza, N. Dobigeon, M. Parente, Q. Du, P. Gader, and J.  Chanussot,
\newblock ``Hyperspectral Unmixing Overview: Geometrical, Statistical and Sparse Regression-Based Approaches'',
\newblock  {\em IEEE Journal of Selected Topics in Applied Earth Observations and Remote Sensing}
 accepted for publication, 2012.


%



\bibitem{BrucksteinElad}
A.~Bruckstein, M.~Elad, and M.~Zibulevsky, ``A non-negative and sparse enough solution of an underdetermined linear system of equations is unique",
{\em IEEE Trans. Inf. Theo.}, vol.~54, pp.~4813--4820, 2008.


\bibitem{CHANG06}
C.-I. Chang, C.-C.~Wu, W.~Liu, and Y.-C.~Ouyang,
\newblock ``A new growing method for simplex-based endmember extraction
  algorithm,''
\newblock {\em IEEE Trans. Geosc. Remote Sensing}, vol. 44, pp. 2804-- 2819, 2006.



\bibitem{chenDonohoSaunders:95}
S. Chen and D. Donoho, and and M. Saunders,
\newblock ``Atomic decomposition by basis pursuit,''
\newblock {\em SIAM review}, vol. 43, no.1, pp. 129--159, 1995.

%
%

\bibitem{Dobigeon:ieeesp:08}
N.~Dobigeon, J.-Y. Tourneret, and C.-I Chang,
\newblock ``Semi-supervised linear spectral unmixing using a hierarchical
  {B}ayesian model for hyperspectral imagery,''
\newblock {\em IEEE Trans. Signal Proc.}, vol. 56, pp.~2684--2695, 2008.


\bibitem{EcksteinBertsekas}
J.~Eckstein, D.~Bertsekas, ``On the Douglas-–Rachford splitting
method and the proximal point algorithm for maximal monotone operators", {\it  Math. Progr.,}
vol.~5, pp.~293–-318, 1992


\bibitem{Gabay}
D.~Gabay and B.~Mercier, ``A dual algorithm for the solution of nonlinear
variational problems via finite-element approximations", {\it Comp. and Math.
 Appl.}, vol.~2, pp.~17--40, 1976.

\bibitem{Glowinski}
R.~Glowinski and A.~Marroco, ``Sur l'approximation, par elements finis d'ordre un, et la
resolution, par penalisation-dualit\'{e} d'une classe de problemes de Dirichlet non lineares",
{\it Rev. Fran\c{c}aise d'Automatique, Inform. Rech. Op\'{e}rationelle},  vol.~9, pp.~41--76, 1975.


\bibitem{conf:daniel:whispers:10}
M.-D.~Iordache, A. Plaza, and J.~Bioucas-Dias,
\newblock "On the Use of Spectral Libraries to perform Sparse Unmixing of Hyperspectral Data",
\newblock in  {\em 2nd IEEE GRSS Workshop on Hyperspectral Image and Signal Processing-WHISPERS'2010}, Reykjvik, Iceland, 2010.


\bibitem{keshava02}
N.~Keshava and J.~F. Mustard, ``Spectral unmixing,'' \emph{IEEE Signal
  Processing Magazine}, vol.~19, no.~1, pp. 44--57, 2002.


\bibitem{Lopes}
M.~Lopes, J.-C.~Wolff, J.~Bioucas-Dias, M.~Figueiredo, ``NIR hyperspectral
unmixing based on a minimum volume criterion for fast and accurate
chemical characterisation of counterfeit tablets", \emph{Analytical Chemistry},
vol. 82, pp. 1462-1469, 2010.

\bibitem{ke00}
D.G.~Manolakis N.~Keshava, J.P.~kerekes and G.A. Shaw,
\newblock ``Algorithm taxonomy for hyperspectral unmixing,''
\newblock {\em Proc. SPIE Vol.4049, Algorithms for Multispectral,
  Hyperspectral, and Ultraspectral Imagery}, vol. VI, pp. 42, 2000.


\bibitem{LM07}
L.~Miao and H.~Qi,
\newblock ``Endmember extraction from highly mixed data using minimum volume
  constrained nonegative matrix factorization,''
\newblock {\em IEEE Trans. Geosc. Remote Sensing}, vol. 45,
  pp. 765--777, 2007.


\bibitem{Moussaoui:nc:08}
S.~Moussaoui, H.~Hauksd\'ottir, F.~Schmidt, C.~Jutten, J.~Chanussot, D.~Brie,
  S.~Dout\'e, and J.~A. Benediksson,
\newblock ``On the decomposition of {M}ars hyperspectral data by {ICA} and
  {B}ayesian positive source separation,''
\newblock {\em Neurocomputing}, 2008,
\newblock accepted.


\bibitem{JNBD05}
J.~Nascimento and J.~Bioucas-Dias,
\newblock ``Does independent component analysis play a role in unmixing hyperspectral data?''
\newblock {\em IEEE Trans. on Geoscience and Remote Sensing}, vol. 43,
pp. 175--187, 2005.

\bibitem{Dias}
J.~Nascimento and J.~Bioucas-Dias,
\newblock ``Hyerspectral unmixing algorithm via dependent component analysis,''
\newblock {\em IEEE Intern. Geoscience and Remote Sensing Symposium}, pp.
  4033--4036, 2007.

\bibitem{plaza04}
R.~Perez A.~Plaza, P.~Martinez and J.~Plaza,
\newblock ``A quantitative and comparative analysis of endmembr extraction
  algorithms from hyperspectral data,''
\newblock {\em IEEE Trans. on Geoscience and Remote Sensing}, vol. 42,
  pp. 650--663, 2004.

\bibitem{plaza02}
R.~Perez A.~Plaza, P.~Martinez and J.~Plaza,
\newblock ``Spatial/spectral endmember extraction by multidimensional
  morphological operations,''
\newblock {\em IEEE Trans. Geosc. Remote Sensing}, vol. 40,
  pp. 2025--2041, 2002.

\bibitem{Roc70} R.~T.~Rockafellar, {\em Convex Analysis}, Princeton
  University Press, Princeton, NJ, 1970.


\bibitem{zare08}
A. Zare and and P. Gader,
\newblock "Hyperspectral band selection and endmember detection using sparsity promoting priors",
\newblock {\em IEEE Geoscience and Remote Sensing Letters}, vol 5., no. 2, pp. 256--260, 2008.



\bibitem{conf:szlam2010split}
A. Szlam, Z.  Guo, Z. and S.  Osher, ``A split Bregman method for non-negative sparsity penalized
least squares with applications to hyperspectral demixing'', {\em IEEE 17th International Conference on
Image Processing (ICIP)}, pp. 1917--1920, 2010.


\bibitem{esser09}
E. Esser,
\newblock ``Applications of lagrangian-based alternating direction methods and connections to
split Bregman,'' Tech. Rep. TR09-31, UCLA CAM, 2009.
\end{thebibliography}
\end{document}